\documentclass[a4paper,12pt]{article}
\usepackage{amsmath,amssymb} \usepackage{euler,eucal}
\usepackage[utf8]{inputenc}
\usepackage[T1]{fontenc}
\usepackage[english]{babel}
\usepackage{tgpagella}
\usepackage[unicode=true,plainpages=false]{hyperref}
\usepackage{graphicx}
\hypersetup{colorlinks=false,pdfborder={0 0 0.1},linkcolor=magenta,anchorcolor=magenta,urlcolor=blue,citecolor=blue,
          pdftitle={Computation of extreme eigenvalues in higher dimensions using block tensor train format},
          pdfauthor={S. V. Dolgov, B. N. Khoromskij, I. V. Oseledets, D. V. Savostyanov},
          pdfkeywords={high-dimensional problems, DMRG, MPS, tensor train format, low-lying eigenstates}
          }
\newcommand*{\doi}[1]{doi: \href{http://dx.doi.org/#1}{#1}}
\usepackage[left=30mm,right=30mm,top=20mm,bottom=30mm,a4paper]{geometry}
\usepackage{algorithmic}
\usepackage[ruled]{algorithm}
\usepackage{amsthm}
\theoremstyle{definition}

\newtheorem{definition}{Definition}
\usepackage{tikz} \usetikzlibrary{positioning,shapes}
\usepackage{pgfplots} \usetikzlibrary{plotmarks}
\usetikzlibrary{calc}
\def\TIKZcoresize{14mm}
\def\TIKZlspan{8mm}
\def\TIKZhspan{7mm}
\input ./tikz/ttdef
\usepackage{bm}
\let\eps=\varepsilon

\let\leq=\leqslant
\let\geq=\geqslant

\def\O{\mathcal{O}}
\def\X{\mathcal{X}}

\def\S{\mathbf{S}}

\def\trace{\mathop{\mathrm{trace}}\nolimits}

\def\vec{\mathop{\mathrm{vec}}\nolimits}
\def\trans{*}
\def\new{\star}
\def\work{\mathop{\mathrm{work}}\nolimits}
\def\time{\mathop{\mathrm{time}}\nolimits}
\def\Span{\mathop{\mathrm{span}}\nolimits}
\def\matvec{\mathop{\mathrm{matvec}}\nolimits}
\begin{document}
\author{
  S. V. Dolgov\footnotemark[2],
  B. N. Khoromskij\footnotemark[2],
  I. V. Oseledets\footnotemark[3],
  D. V. Savostyanov\footnotemark[3]~\footnotemark[4]
}
\title{Computation of extreme eigenvalues in higher dimensions using block tensor train format%
 \thanks{Partially supported by
         RFBR grants 11-01-00549-a, 12-01-33013, 12-01-00546-a,
         12-01-91333-nnio-a,
         Rus. Fed. Gov. project 16.740.12.0727
         at the Institute of Numerical Mathematics RAS
         and EPSRC grant EP/H003789/1 at the University of Southampton.
         }}
\date{June 09, 2013}
\maketitle
\renewcommand{\thefootnote}{\fnsymbol{footnote}}
\footnotetext[2]{Max-Planck Institute for Mathematics in the Sciences, Inselstrasse 22, Leipzig 04103, Germany  ({\tt bokh@mis.mpg.de, dolgov@mis.mpg.de})}
\footnotetext[3]{Institute of Numerical Mathematics of Russian Academy of Sciences, Gubkina 8, Moscow 119333, Russia ({\tt ivan.oseledets@gmail.com})}
\footnotetext[4]{University of Southampton, School of Chemistry, Highfield Campus, Southampton SO17 1BJ, United Kingdom  ({\tt dmitry.savostyanov@gmail.com})}
\renewcommand{\thefootnote}{\arabic{footnote}}

\begin{abstract}
We consider an approximate computation of several minimal eigenpairs of large Hermitian matrices which come from high--dimensional problems.
We use the tensor train format (TT) for vectors and matrices to overcome the curse of dimensionality and make storage and computational cost feasible.
Applying a block version of the TT format to several vectors simultaneously, we compute the low--lying eigenstates of a system by minimization of a block Rayleigh quotient performed in an alternating fashion for all dimensions.
For several numerical examples, we compare the proposed method with the deflation approach when the low--lying eigenstates are computed one-by-one, and also with the variational algorithms used in quantum physics.
\end{abstract}

\section{Introduction}

High-dimensional problems are notoriously difficult to solve by standard numerical techniques due to the \emph{curse of dimensionality} --- the complexity grows exponentially with the number of degrees of freedom.
The problems of such kind arise in many different applications in physics, chemistry, biology and engineering, but their study in numerical linear algebra has begun quite recently.

There are not many techniques capable of solving high--dimensional problems efficiently.
The most prominent among them are Monte Carlo and quasi Monte Carlo methods, best N-term approximations,  and advanced discretization methods such as sparse grids and radial basis functions.
However, all of these methods have their own disadvantages.
For example, it is difficult to achieve high accuracy using the Monte Carlo approach, and sparse grid techniques require sophisticated analytical and algebraic manipulations and still suffer (in a milder way though) from the curse of dimensionality, which make them inapplicable for $d \gtrsim 10.$

One of the most fruitful ideas for solving high-dimensional problems is the idea of \emph{separation of variables}.
For two variables it boils down to the celebrated Schmidt decomposition, which is known on a discrete level as the singular value decomposition (SVD), a particular low--rank decomposition of a matrix.
Different generalizations of this idea to higher dimensions, most notable are the canonical (CP) and Tucker formats, have been studied motivated by the applications in data analysis, particularly chemometrics, see the review~\cite{kolda-review-2009}.
These classical formats have their drawbacks as well, e.g. the canonical format is in general not stable to perturbations, and the Tucker format suffers from the curse of dimensionality.
Nevertheless, in many applications the canonical representation can be computed efficiently using, e.g.~\emph{greedy algorithms}~\cite{chinesta-greedy-2006,maday-greedy-2009,Temlyakov-greedy-2011} or by a multigrid accelerated reduced higher order
SVD combined with the Tucker format \cite{khor-ml-2009}, and for the Tucker format a quasioptimal approximation can be computed using the SVD algorithm~\cite{lathauwer-svd-2000}.

Efficient methods for quantum many-body systems are based on low--parametric tensor product formats. One of the most successful approaches,
the \emph{density matrix renormalization group} (DMRG) \cite{white-dmrg-1992,white-dmrg-1993} is an optimization technique that
uses the matrix product state (MPS) representation, see the review \cite{schollwock-review-2011}.
The MPS and DMRG are described in a problem--specific language, and despite they became the methods of choice for many applications
in the solid state physics and quantum chemistry,  they were unknown in numerical analysis.

Looking for more efficient dimensionality reduction schemes, two groups in the numerical linear algebra community have re-discovered independently successful tensor formats under different names, namely the Tree-Tucker~\cite{ot-tt-2009}, tensor train (TT)~\cite{osel-tt-2011} and hierarchical Tucker (HT)~\cite{hk-ht-2009,gras-hsvd-2010} formats.
The equivalence of the TT and MPS format has been shortly discovered and reported in~\cite{holtz-ALS-DMRG-2012}.
This connection is very beneficial and fruitful: the idea of the DMRG algorithm has been applied to different kinds of problem in numerical analysis: approximate solution of linear systems \cite{holtz-ALS-DMRG-2012, DoOs-dmrg-solve-2011}, solution of eigenvalue problems \cite{khos-dmrg-2010}, dynamics~\cite{DKhOs-parabolic1-2012}, cross interpolation \cite{so-dmrgi-2011proc}.
At the same time, new tensor formats have been proposed, e.g. the \emph{quantized} tensor train (QTT)~\cite{khor-qtt-2011,osel-2d2d-2010}, and the QTT-Tucker~\cite{dk-qtt-tucker-2013}, and several new results have been obtained for eigenproblems~\cite{lebedeva-tensornd-2011}, solution of linear systems \cite{ds-amr1-2013, ds-amr2-2013}, multidimensional convolution~\cite{khor-acc-2010,hackbusch-conv-2011,khkaz-conv-2011}, multidimensional Fourier transform~\cite{dks-ttfft-2012},  interpolation~\cite{sav-qott-2013pre}.

The use of different tensor product formats for quantum chemistry problems has been considered in~\cite{khor-tuckertype-2007,fkst-chem-2008,khor-ml-2009,khst-eigen-2012,khor-chem-2011,vkhs-2el-2013}.
This paper can be considered as a natural extension of the work~\cite{khos-dmrg-2010}, where the computation of a single extreme eigenvalue using the (quantized) tensor train format was revisited.
Instead of computing one eigenvalue we want to compute approximations to several extreme eigenpairs simultaneously.
Such problems appear, for example, in the computation of excited states in physics and chemistry.
This generalization is not straightforward, taking into account the special structure of the manifold where the solution is sought.
Our approach consists in two key components.
First, all eigenstates are represented in the TT-format in a simultaneous manner (so-called \emph{block TT-format}), and the minimization problem is formulated in a rigorous way as a minimization of the Block Rayleigh quotient.
The optimization problem is then solved by means of alternating optimization with a small but very important trick: instead of using one-and-the same representation for all steps, the auxiliary index corresponding to the eigenstate number is always associated with the currently optimized TT core.
The local problem at each iteration step becomes a linear block eigenvalue problem which can be treated by standard iterative techniques.

At the moment the draft of this paper was written, Frank Verstraete and Iztok Pizorn kindly informed us about the work~\cite{pizorn-eigb-2012} in the MPS community,  which also addresses solution of the block Rayleigh minimization problem in the alternating framework, but the methods described in \cite{pizorn-eigb-2012} differ from the algorithm proposed in this paper.

Our algorithm has a better asymptotic complexity with respect to the mode size, but adapts the TT ranks during the iterations like in the DMRG scheme.
This can be particularly important for solving high--dimensional problems with large mode sizes.

The paper is organized as follows. In Section \ref{sec:def} definitions of the tensor train format are introduced.
In Section \ref{sec:alg} the algorithm is presented, and the complexity is analyzed and compared to the algorithms used in quantum physics.
Sections~\ref{sec:lap},~\ref{sec:hh},~\ref{sec:spin} contain numerical experiments for the particle in a box, the Henon--Heiles potential, and the Heisenberg model,
including the comparison of the computational speed with the publicly available DMRG implementations.

\section{Notation, definitions and preliminaries} \label{sec:def}
We consider the eigenproblem $AX=X\Lambda,$ with the Hermitian matrix $A=A^\trans,$ and we are interested in $B$ extreme eigenvalues $\lambda_b$ and their eigenvectors $x_b,$ for $b=1,\ldots,B.$
This problem is equivalent to the minimization of the block Rayleigh quotient
\begin{equation}\label{eq:rq}
 \trace(X^{\trans} A X) \to \min, \qquad \mbox{s.t.}\quad X^{\trans} X = I,
\end{equation}
where $X=[x_b]_{b=1}^B$ contains the orthogonal eigenvectors.

We assume that the problem has a tensor-product structure, i.e. all eigenvectors involved can be associated with $d$-dimensional tensors.
Specifically, the elements of a vector $x=[x(i)]_{i=1}^N$ can be enumerated with $d$ \emph{mode indices} $i_1,\ldots,i_d$ by a linear map $i=\overline{i_1\ldots i_d}.$
The mode indices run through $i_k=1,\ldots,n_k,$ where $n_k$ are referred to as the~\emph{mode sizes} for $k=1,\ldots,d.$
Naturally, $N=n_1\ldots n_d$, and if all mode sizes are of the same order $n_k\sim n,$ the number of unknowns grows exponentially with the dimension, $N\sim n^d.$
To make the problem tractable, we use the~\emph{tensor train} (TT) format~\cite{osel-tt-2011}, defined as follows,
\begin{equation}\label{eq:tt}
 \begin{split}
 x(i) &  = x(\overline{i_1\ldots i_d})
       = X^{(1)}(i_1) \ldots  X^{(d)}(i_d)
 \\   & = \sum_{\bm{\alpha}}  X^{(1)}_{\alpha_1}(i_1) \ldots  X^{(k-1)}_{\alpha_{k-2},\alpha_{k-1}}(i_{k-1})X^{(k)}_{\alpha_{k-1},\alpha_k}(i_k)X^{(k+1)}_{\alpha_{k},\alpha_{k+1}}(i_{k+1})\ldots X^{(d)}_{\alpha_{d-1}}(i_d).
 \end{split}
\end{equation}
\begin{figure}[t]
 \begin{center}
  \resizebox{.95\textwidth}{!}{ \def\ss{\scriptsize}
 \def\st{\scriptstyle}
 \def\L{\TIKZlspan}
 \def\H{\TIKZhspan}
 \begin{tikzpicture}[x=10mm,y=-10mm]
   \node[core] (a5)                   {$X^{(k)}$};    
   \node[core] (a6) [right=\L of a5]  {$X^{(k+1)}$};
   \node[void] (a7) [right=\L of a6]  {};
   \node[core] (a8) [right=\L of a7]  {$X^{(d)}$};
   \node[core] (a4) [left =\L of a5]  {$X^{(k-1)}$};
   \node[void] (a3) [left =\L of a4]  {};
   \node[core] (a2) [left =\L of a3]  {$X^{(2)}$};
   \node[core] (a1) [left= \L of a2]  {$X^{(1)}$};
   
   \path[bond]        (a1)      -- node[above,midway]  {${\st\alpha_1}$}       (a2);
   \path[bond]        (a2)      -- node[above,midway]  {${\st\alpha_2}$}       (a3.west);
   \path[bond,dashed] (a3.west) --                                             (a3.east);
   \path[bond]        (a3.east) -- node[above,midway]  {${\st\alpha_{k-2}}$}   (a4);
   \path[bond]        (a4)      -- node[above,midway]  {${\st\alpha_{k-1}}$}   (a5);
   \path[bond]        (a5)      -- node[above,midway]  {${\st\alpha_{k}}$}     (a6);
   \path[bond]        (a6)      -- node[above,midway]  {${\st\alpha_{k+1}}$}   (a7.west);
   \path[bond,dashed] (a7.west) --                                             (a7.east);
   \path[bond]        (a7.east) -- node[above,midway]  {${\st\alpha_d}$}       (a8);
   
   \path[bond]        (a1.south)   -- node[left,midway] {${\st i_1}$}     +(0,+0.5);
   \path[bond]        (a2.south)   -- node[left,midway] {${\st i_2}$}     +(0,+0.5);
   \path[bond]        (a4.south)   -- node[left,midway] {${\st i_{k-1}}$} +(0,+0.5);
   \path[bond]        (a5.south)   -- node[left,midway] {${\st i_{k}}$}   +(0,+0.5);
   \path[bond]        (a6.south)   -- node[left,midway] {${\st i_{k+1}}$} +(0,+0.5);
   \path[bond]        (a8.south)   -- node[left,midway] {${\st i_d}$}     +(0,+0.5);
   
 \end{tikzpicture}
 }
  \vskip 5mm
  \resizebox{.95\textwidth}{!}{ \def\ss{\scriptsize}
 \def\st{\scriptstyle}
 \def\L{\TIKZlspan}
 \def\H{\TIKZhspan}
 \begin{tikzpicture}[x=10mm,y=-10mm]
   \node[core] (a5)                   {$A^{(k)}$};    
   \node[core] (a6) [right=\L of a5]  {$A^{(k+1)}$};
   \node[void] (a7) [right=\L of a6]  {};
   \node[core] (a8) [right=\L of a7]  {$A^{(d)}$};
   \node[core] (a4) [left =\L of a5]  {$A^{(k-1)}$};
   \node[void] (a3) [left =\L of a4]  {};
   \node[core] (a2) [left =\L of a3]  {$A^{(2)}$};
   \node[core] (a1) [left= \L of a2]  {$A^{(1)}$};
   
   \path[bond]        (a1)      -- node[above,midway]  {${\st\alpha_1}$}       (a2);
   \path[bond]        (a2)      -- node[above,midway]  {${\st\alpha_2}$}       (a3.west);
   \path[bond,dashed] (a3.west) --                                             (a3.east);
   \path[bond]        (a3.east) -- node[above,midway]  {${\st\alpha_{k-2}}$}   (a4);
   \path[bond]        (a4)      -- node[above,midway]  {${\st\alpha_{k-1}}$}   (a5);
   \path[bond]        (a5)      -- node[above,midway]  {${\st\alpha_{k}}$}     (a6);
   \path[bond]        (a6)      -- node[above,midway]  {${\st\alpha_{k+1}}$}   (a7.west);
   \path[bond,dashed] (a7.west) --                                             (a7.east);
   \path[bond]        (a7.east) -- node[above,midway]  {${\st\alpha_d}$}       (a8);
   
   \path[bond]        (a1.north)   -- node[left,midway] {${\st i_1}$}         +(0,-0.5);
   \path[bond]        (a2.north)   -- node[left,midway] {${\st i_2}$}         +(0,-0.5);
   \path[bond]        (a4.north)   -- node[left,midway] {${\st i_{k-1}}$}     +(0,-0.5);
   \path[bond]        (a5.north)   -- node[left,midway] {${\st i_k}$}         +(0,-0.5);
   \path[bond]        (a6.north)   -- node[left,midway] {${\st i_{k+1}}$}     +(0,-0.5);
   \path[bond]        (a8.north)   -- node[left,midway] {${\st i_d}$}         +(0,-0.5);
   
   \path[bond]        (a1.south)   -- node[right,midway] {${\st j_1}$}        +(0,+0.5);
   \path[bond]        (a2.south)   -- node[right,midway] {${\st j_2}$}        +(0,+0.5);
   \path[bond]        (a4.south)   -- node[right,midway] {${\st j_{k-1}}$}    +(0,+0.5);
   \path[bond]        (a5.south)   -- node[right,midway] {${\st j_k}$}        +(0,+0.5);
   \path[bond]        (a6.south)   -- node[right,midway] {${\st j_{k+1}}$}    +(0,+0.5);
   \path[bond]        (a8.south)   -- node[right,midway] {${\st j_d}$}        +(0,+0.5);
   
 \end{tikzpicture}
 }
 \end{center}
 \caption{Above: tensor train format~\eqref{eq:tt} shown as a linear tensor network.
          Below: tensor train format~\eqref{eq:att}  for a matrix in higher dimensions.
          The boxes are tensors with lines (legs) denoting indices.
          Each bond between two tensors assumes a summation over the joint index.}
 \label{fig:tt}
\end{figure}
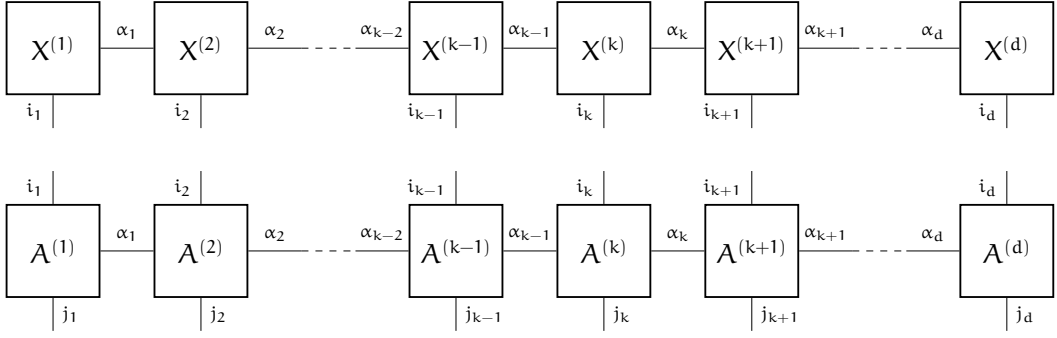
Here and later we write equations in the \emph{elementwise notation}, i.e. assume they hold for all possible values of all free indices.
The summation over $\bm{\alpha}=(\alpha_1,\ldots,\alpha_{d-1})$ assumes the summation over all possible values of all \emph{auxiliary} (or bond) indices $\alpha_k=1,\ldots,r_k,$ where numbers $r_1,\ldots,r_{d-1}$ are referred to as the \emph{tensor train ranks} (TT--ranks).
Each $X^{(k)}(i_k)$ is an $r_{k-1}\times r_k$ matrix,  i.e. each entry of a vector $x=[x(i)]$ is represented by a product of $d$ matrices in the right--hand side.
The three--dimensional arrays $X^{(k)}=[X^{(k)}_{\alpha_{k-1},\alpha_k}(i_k)]$  of size $r_{k-1} \times n_k \times r_k$ are referred to as the \emph{TT--cores}.
In the representation~\eqref{eq:tt}, the vector $x$ is seen as a~\emph{vectorization} of a $d$--tensor, which explains the name of the format.
The tensor train format is a \emph{linear tensor network}, and can be illustrated as a graph, see Fig.~\ref{fig:tt}.

In this paper we represent all eigenvectors of interest~\emph{simultaneously} by the \emph{block} tensor train format.
\begin{definition}[block TT-format]
The vectors $x_b(i)$ are said to be in the block--TT format, if
\begin{equation}\label{eq:btt}
 \begin{split}
 x_b(i) & = x_b(\overline{i_1\ldots i_d})
        = X^{(1)}(i_1) \ldots X^{(p-1)}(i_{p-1}) \hat X^{(p)}(i_p,b) X^{(p+1)}(i_{p+1}) \ldots X^{(d)}(i_d).
 \end{split}
\end{equation}
\end{definition}
The choice of the mode $p$ which \emph{carries} the index $b$ is not fixed --- we will move it back and forth during the optimization.
When the position $p$ is chosen, it means that the matrix $\hat X^{(p)}(i_p,b)$ is additionally parametrized by the index $b.$
The `block' TT--core $\hat X^{(p)}=[\hat X^{(p)}_{\alpha_{p-1},\alpha_p}(i_p,b)]$ is now a tensor with four indices.

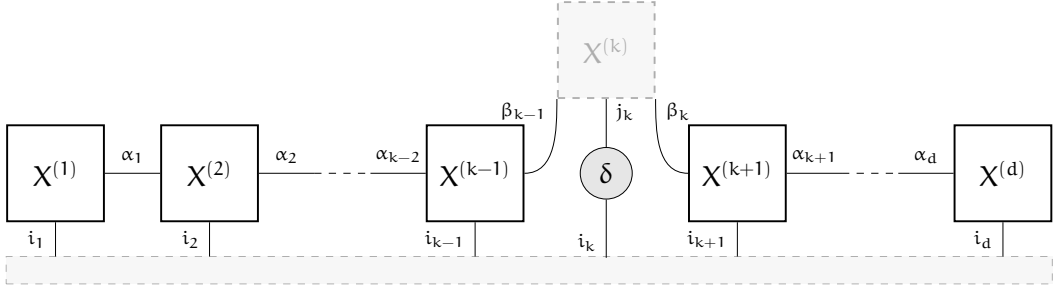
\begin{figure}[t]
 \begin{center}
  \resizebox{.95\textwidth}{!}{ \def\ss{\scriptsize}
 \def\st{\scriptstyle}
 \def\L{\TIKZlspan}
 \def\H{\TIKZhspan}
 \begin{tikzpicture}[x=10mm,y=-10mm]
   \node[shape=circle,minimum size=1em,draw=black,fill=black!10]  (a5) {$\delta$};
   \node[core,dashed,fill=black!10,opacity=.3] (b5) [above=\H of a5] {{$X^{(k)}$}};
   \node[core] (a6) [right=\L of a5]  {$X^{(k+1)}$};
   \node[void] (a7) [right=\L of a6]  {};
   \node[core] (a8) [right=\L of a7]  {$X^{(d)}$};
   \node[core] (a4) [left =\L of a5]  {$X^{(k-1)}$};
   \node[void] (a3) [left =\L of a4]  {};
   \node[core] (a2) [left =\L of a3]  {$X^{(2)}$};
   \node[core] (a1) [left= \L of a2]  {$X^{(1)}$};

   \filldraw[dashed,fill=black!10,opacity=0.3] ($(a1.south west)+(0,+0.5)$) rectangle ($(a8.south east)+(0,+.9)$);

   \path[bond]        (a1)      -- node[above,midway]  {${\st\alpha_1}$}       (a2);
   \path[bond]        (a2)      -- node[above,midway]  {${\st\alpha_2}$}       (a3);
   \path[bond,dashed] (a3.west) --                                             (a3.east);
   \path[bond]        (a3.east) -- node[above,midway]  {${\st\alpha_{k-2}}$}   (a4);
   \draw[]            (a4.east)  to[out=right,in=down] node[left,very near end] {${\st\beta_{k-1}}$} (b5.south west);
   \draw[]            (a6.west)  to[out=left,in=down]  node[right,very near end] {${\st\beta_{k}}$}   (b5.south east);
   \path[bond]        (a6)      -- node[above,midway]  {${\st\alpha_{k+1}}$}   (a7.west);
   \path[bond,dashed] (a7.west) --                                             (a7.east);
   \path[bond]        (a7.east) -- node[above,midway]  {${\st\alpha_d}$}       (a8);

   \path[bond]        (a1.south)   -- node[left,midway] {${\st i_1}$}         +(0,+0.5);
   \path[bond]        (a2.south)   -- node[left,midway] {${\st i_2}$}         +(0,+0.5);
   \path[bond]        (a4.south)   -- node[left,midway] {${\st i_{k-1}}$}     +(0,+0.5);
   \draw[]            (a5.north) to[out=up,in=down]  node[right,near end] {${\st j_k}$}  (b5.south);
   \draw[]            (a5.south) to[out=down,in=up]  node[left,near end] {${\st i_k}$}  ($(a5.south)+(0,+.85)$);
   \path[bond]        (a6.south)   -- node[left,midway] {${\st i_{k+1}}$}     +(0,+0.5);
   \path[bond]        (a8.south)   -- node[left,midway] {${\st i_d}$}         +(0,+0.5);
 \end{tikzpicture}
 }
 \end{center}
 \caption{The frame matrix~\eqref{eq:frame} maps a TT core (above) into a large vector (below).}
 \label{fig:frame}
\end{figure}

Following~\cite{ushmaev-tt-2013}, we define the \emph{interfaces} $X^{<k}$ of size ${n_1\ldots n_{k-1} \times r_{k-1}}$ and $X^{>k}$ of size ${r_k \times n_{k+1}\ldots n_d}$ as follows
\begin{equation}\label{eq:iface}
 \begin{split}
     X^{<k}(\overline{i_1 i_2 \ldots i_{k-1}}, \beta_{k-1}) &
    = \sum_{\alpha_1 \ldots \alpha_{k-2}}
       X^{(1)}_{\alpha_1}(i_1) X^{(2)}_{\alpha_1\alpha_2}(i_2) \ldots X^{(k-1)}_{\alpha_{k-2},\beta_{k-1}}(i_{k-1}), \\
  X^{> k}(\beta_k,\overline{i_{k+1}\ldots i_{d-1}i_d}) &
  = \sum_{\alpha_{k+1} \ldots \alpha_{d-1}}
     X^{(k+1)}_{\beta_k,\alpha_{k+1}}(i_{k+1}) \ldots X^{(d-1)}_{\alpha_{d-2},\alpha_{d-1}}(i_{d-1}) X^{(d)}_{\alpha_{d-1}}(i_d).
 \end{split}
\end{equation}
Using the interfaces we introduce the~\emph{frame matrix} (see Fig.~\ref{fig:frame}) as follows
\begin{equation}\label{eq:frame}
 \begin{split}
 \X_{\neq k} & =  X^{<k} \otimes I_{n_k} \otimes \left(X^{>k} \right)^\top, \\
 \X_{\neq k}&(\overline{i_1\ldots i_d},  \overline{\beta_{k-1} j_k \beta_k})
    \\ & =  \sum_{ \substack{ \alpha_1,\ldots,\alpha_{k-2} \\ \alpha_{k+1} \ldots \alpha_d }}
    X^{(1)}_{\alpha_1}(i_1) \ldots X^{(k-1)}_{\alpha_{k-2},\beta_{k-1}}(i_{k-1}) \delta(i_k,j_k)X^{(k+1)}_{\beta_k,\alpha_{k+1}}(i_{k+1}) \ldots X^{(d)}_{\alpha_{d-1}}(i_d).
 \end{split}
\end{equation}

The tensor train format~\eqref{eq:tt} is polylinear with respect to  the TT--cores, and in particular the block-TT format~\eqref{eq:btt} is linear w.r.t. the block core $\hat X^{(p)}.$
We can write it using the frame matrix as follows,
\begin{equation}\label{eq:lin}
x = \X_{\neq k}\hat x^{(k)}, \qquad X = \X_{\neq p}\hat X^{(p)},
\end{equation}
where the vector  $x^{(k)}$ of size $r_{k-1}n_kr_k$ is a vectorization of the TT--core $X^{(k)}$,
$$
x^{(k)} = \vec X^{(k)}, \qquad x^{(k)}(\overline{\beta_{k-1}j_k\beta_k}) = X^{(k)}_{\beta_{k-1},\beta_k}(j_k),
$$
and the block TT--core $\hat X^{(p)}=[\hat x^{(p)}_b]$ is seen as a $r_{p-1}n_pr_p \times B$ matrix, which encapsulates vectorizations $\hat x^{(p)}_b=\vec\hat X^{(p)}(b)$ as columns.

The TT--core $\hat X^{(p)}$ can be also reshaped into a $r_{p-1}n_p \times B r_p$ matrix by arranging the indices in the form $\hat X^{(p)} = [\hat X^{(p)}(\overline{\alpha_{p-1}i_p},\overline{b\alpha_{p}})],$ where the overline shows that the indices $\alpha_{p-1}$ and $i_p$ are combined in a single multi-index, as well as $b$ and $\alpha_p.$
We can write the low--rank decomposition of this matrix
\begin{equation}\label{eq:lr}
 \hat X^{(p)}_{\alpha_{p-1},\alpha_p}(i_p,b) = \hat X^{(p)}(\overline{\alpha_{p-1}i_p},\overline{b\alpha_{p}}) = \sum_{\alpha'_p=1}^{r'_p} X^{(p)}_{\alpha_{p-1},\alpha'_p}(i_p) G_{\alpha'_p,\alpha_p}(b),
\end{equation}
where the new matrix $G(b)$ of size $r'_p\times r_p$ appears, which carries the index $b.$
Substituting~\eqref{eq:lr}  into~\eqref{eq:btt} as follows,
\begin{equation}\label{eq:gtt}
 x_b(i) = x_b(\overline{i_1\ldots i_d})
        = X^{(1)}(i_1) \ldots  X^{(p)}(i_p) G(b) X^{(p+1)}(i_{p+1}) \ldots X^{(d)}(i_d),
\end{equation}
we obtain the regular TT/MPS format for $X=[x_b(i)],$ where the index $i$ is defined by $d$ mode indices $i_1,\ldots,i_d$ and $b$ is considered as the additional index which is placed between $i_p$ and $i_{p+1}.$
This format with $G(b)$ in the rightmost position has been applied for the eigenproblem solution in~\cite{lebedeva-tensornd-2011}, and is also considered in~\cite{schollwock-review-2011}.

Using the decomposition~\eqref{eq:lr}, we can move the index $b$ back and forth the tensor train.
Indeed, if we multiply $G(b)$ by $X^{(p+1)}$ and reshape the result, we obtain $\hat X^{(p+1)}_b$ and recover~\eqref{eq:btt} with the `block' TT--core at the position $p+1.$
We see that after this multiplication the TT--rank $r_p$ is replaced by $r_p'$ and in general $r_p'\neq r_p.$
This means that the TT--ranks depend on the position $p$ of the TT--core which carries the index $b.$
The natural bound for the considered example is $r'_p\leq r_p B,$ which technically allows the ranks to grow exponentially during the iterations.
For the optimization algorithm to be efficient, special measures should be taken to keep the TT--ranks moderate.

To reduce (truncate) the TT--rank, we can compute the low--rank decomposition~\eqref{eq:lr} approximately.
The approximation of the lowest rank within the prescribed accuracy level can be computed by the \emph{singular value decomposition} (SVD).
The perturbation introduced to $\hat X^{(p)}$ results in the perturbation to the whole matrix of eigenstates $X,$ which can be amplified by the norm of the frame matrix $\X_{\neq p}$ due to the linearity~\eqref{eq:lin}.
If $\X_{\neq p}^\trans \X_{\neq p}=I,$ the Frobenius norms of the local and the global perturbations are the same.
For a given tensor the orthogonality of the frame matrix can be achieved using the non--uniqueness of the TT format, by implying the \emph{left--orthogonality} (left--normalization) constrains to $X^{(1)},\ldots,X^{(p-1)}$ and the \emph{right--orthogonality} (right--normalization) to $X^{(p+1)},\ldots,X^{(d)}.$
\begin{definition}[TT--orthogonality]
    The TT--core $X^{(k)}$ is called \emph{left--} or \emph{right--orthogonal}, if, respectively,
    $$
    \sum_{i_k} (X^{(k)}(i_k))^\trans X^{(k)}(i_k) = I, \qquad
    \sum_{i_k} X^{(k)}(i_k) (X^{(k)}(i_k))^\trans = I.
    $$
\end{definition}
It can be shown \cite{osel-tt-2011,schollwock-review-2011} that if $X^{(1)},\ldots,X^{(p-1)}$ are left-orthogonal, then the interface matrix $X^{<p}$ defined by~\eqref{eq:iface} has orthonormal columns.
Analogously, right-orthogonality of the TT-cores $X^{(p+1)},\ldots,X^{(d)}$ implies that the rows of $X^{>p}$ are orthonormal.
Together these conditions provide the orthogonality of the frame matrix $\X_{\neq p}$ defined by~\eqref{eq:frame}, which allows the full accuracy control in the rank truncation step.

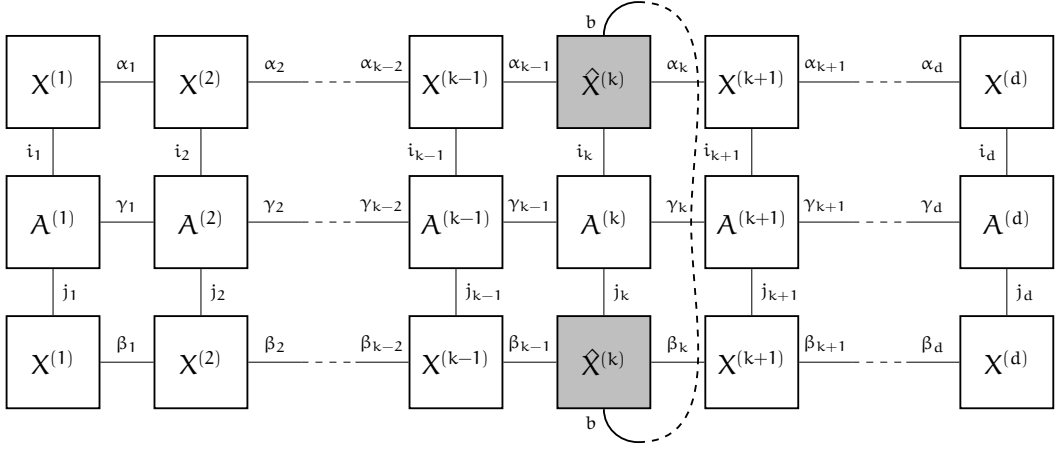
\begin{figure}[t]
 \begin{center}
  \resizebox{.95\textwidth}{!}{ \def\ss{\scriptsize}
 \def\st{\scriptstyle}
 \def\L{\TIKZlspan}
 \def\H{\TIKZhspan}
 \begin{tikzpicture}[x=10mm,y=-10mm]
   \node[core ] (a5)                   {$A^{(k)}$};    
   \node[corex] (x5) [above=\H of a5]  {$\hat X^{(k)}$};
   \node[corex] (y5) [below=\H of a5]  {$\hat X^{(k)}$};
   
   \node[core] (a6) [right=\L of a5]  {$A^{(k+1)}$};
   \node[void] (a7) [right=\L of a6]  {};
   \node[core] (a8) [right=\L of a7]  {$A^{(d)}$};
   \node[core] (a4) [left =\L of a5]  {$A^{(k-1)}$};
   \node[void] (a3) [left =\L of a4]  {};
   \node[core] (a2) [left =\L of a3]  {$A^{(2)}$};
   \node[core] (a1) [left= \L of a2]  {$A^{(1)}$};
   
   \node[core] (x6) [right=\L of x5]  {$X^{(k+1)}$};
   \node[void] (x7) [right=\L of x6]  {};
   \node[core] (x8) [right=\L of x7]  {$X^{(d)}$};
   \node[core] (x4) [left =\L of x5]  {$X^{(k-1)}$};
   \node[void] (x3) [left =\L of x4]  {};
   \node[core] (x2) [left =\L of x3]  {$X^{(2)}$};
   \node[core] (x1) [left= \L of x2]  {$X^{(1)}$};
   
   \node[core] (y6) [right=\L of y5]  {$X^{(k+1)}$};
   \node[void] (y7) [right=\L of y6]  {};
   \node[core] (y8) [right=\L of y7]  {$X^{(d)}$};
   \node[core] (y4) [left =\L of y5]  {$X^{(k-1)}$};
   \node[void] (y3) [left =\L of y4]  {};
   \node[core] (y2) [left =\L of y3]  {$X^{(2)}$};
   \node[core] (y1) [left= \L of y2]  {$X^{(1)}$};

   \path[bond]        (a1) -- node[left,midway]  {${\st i_1}$}     (x1);
   \path[bond]        (a2) -- node[left,midway]  {${\st i_2}$}     (x2);
   \path[bond]        (a4) -- node[left,midway]  {${\st i_{k-1}}$} (x4);
   \path[bond]        (a5) -- node[left,midway]  {${\st i_k}$}     (x5);
   \path[bond]        (a6) -- node[left,midway]  {${\st i_{k+1}}$} (x6);
   \path[bond]        (a8) -- node[left,midway]  {${\st i_d}$}     (x8);
   
   \path[bond]        (a1) -- node[right,midway] {${\st j_1}$}     (y1);
   \path[bond]        (a2) -- node[right,midway] {${\st j_2}$}     (y2);
   \path[bond]        (a4) -- node[right,midway] {${\st j_{k-1}}$} (y4);
   \path[bond]        (a5) -- node[right,midway] {${\st j_k}$}     (y5);
   \path[bond]        (a6) -- node[right,midway] {${\st j_{k+1}}$} (y6);
   \path[bond]        (a8) -- node[right,midway] {${\st j_d}$}     (y8);
   
   \path[bond]        (a1)      -- node[above,midway]  {${\st\gamma_1}$}       (a2);
   \path[bond]        (a2)      -- node[above,midway]  {${\st\gamma_2}$}       (a3.west);
   \path[bond,dashed] (a3.west) --                                             (a3.east);
   \path[bond]        (a3.east) -- node[above,midway]  {${\st\gamma_{k-2}}$}   (a4);
   \path[bond]        (a4)      -- node[above,midway]  {${\st\gamma_{k-1}}$}   (a5);
   \path[bond]        (a5)      -- node[above,midway]  {${\st\gamma_{k}}$}     (a6);
   \path[bond]        (a6)      -- node[above,midway]  {${\st\gamma_{k+1}}$}   (a7.west);
   \path[bond,dashed] (a7.west) --                                             (a7.east);
   \path[bond]        (a7.east) -- node[above,midway]  {${\st\gamma_d}$}       (a8);
   
   \path[bond]        (x1)      -- node[above,midway]  {${\st\alpha_1}$}       (x2);
   \path[bond]        (x2)      -- node[above,midway]  {${\st\alpha_2}$}       (x3.west);
   \path[bond,dashed] (x3.west) --                                             (x3.east);
   \path[bond]        (x3.east) -- node[above,midway]  {${\st\alpha_{k-2}}$}   (x4);
   \path[bond]        (x4)      -- node[above,midway]  {${\st\alpha_{k-1}}$}   (x5);
   \path[bond]        (x5)      -- node[above,midway]  {${\st\alpha_{k}}$}     (x6);
   \path[bond]        (x6)      -- node[above,midway]  {${\st\alpha_{k+1}}$}   (x7.west);
   \path[bond,dashed] (x7.west) --                                             (x7.east);
   \path[bond]        (x7.east) -- node[above,midway]  {${\st\alpha_d}$}       (x8);
   
   \path[bond]        (y1)      -- node[above,midway]  {${\st\beta_1}$}       (y2);
   \path[bond]        (y2)      -- node[above,midway]  {${\st\beta_2}$}       (y3.west);
   \path[bond,dashed] (y3.west) --                                            (y3.east);
   \path[bond]        (y3.east) -- node[above,midway]  {${\st\beta_{k-2}}$}   (y4);
   \path[bond]        (y4)      -- node[above,midway]  {${\st\beta_{k-1}}$}   (y5);
   \path[bond]        (y5)      -- node[above,midway]  {${\st\beta_{k}}$}     (y6);
   \path[bond]        (y6)      -- node[above,midway]  {${\st\beta_{k+1}}$}   (y7.west);
   \path[bond,dashed] (y7.west) --                                            (y7.east);
   \path[bond]        (y7.east) -- node[above,midway]  {${\st\beta_d}$}       (y8);

   \coordinate (xx) at ($(x5.north east)+(-2mm,+5mm)$);
   \coordinate (yy) at ($(y5.south east)+(-2mm,-5mm)$);
   \coordinate (zz) at ($(a6.west)+(left:2mm)$);
   \draw[thick] (x5.north) to[out=up,  in=180] node[left,near start] {$\st b$} (xx); 
   \draw[thick] (y5.south) to[out=down,in=180] node[left,near start] {$\st b$} (yy);
   \draw[thick,dashed] (xx) to[out=0,in= 90] (zz);
   \draw[thick,dashed] (yy) to[out=0,in=-90] (zz);
 \end{tikzpicture}
 }
 \end{center}
 \caption{Tensor network corresponding to the block Rayleigh quotient~\eqref{eq:rq} with the matrix $A$ and vectors $X=[x_b]$ given in the tensor train format.
 }
 \label{fig:xax}
\end{figure}

For a given TT format the TT--orthogonality can be implied constructively by the subsequent orthogonalization of the cores, see~\cite{osel-tt-2011,schollwock-review-2011}.
The procedure is very cheap and requires only $\O(nr^3)$ operations for each TT--core to be normalized.
This operation is never a bottleneck in our algorithms, so we always assume that for a chosen TT--core $\hat X^{(p)}$ the frame matrix is orthogonal, without going into the minor detail.

The TT format for a matrix is written as follows,
\begin{equation}\label{eq:att}
 A(i,j)=A(\overline{i_1\ldots i_d},\overline{j_1 \ldots j_d}) = A^{(1)}(i_1,j_1) \ldots A^{(d)}(i_d,j_d).
\end{equation}
The multiplication of a matrix by a vector in the TT format has been discussed in detail in~\cite{schollwock-dmrg-mps-2010,Os-mvk2-2011}.
The~\emph{tensor network} which represents the quadratic function in~\eqref{eq:rq} is shown in Fig.~\ref{fig:xax}.

\section{Minimization of the block Rayleigh quotient in the block--TT format} \label{sec:alg}
We formulate the eigenproblem via the Rayleigh quotient optimization~\eqref{eq:rq}, which is restrictively large when $d$ increases.
To make the problem tractable, the minimization over the whole space is substituted by the optimization over the manifold of tensors in the TT-format.
The most natural approach is the alternating least squares (ALS)-type algorithm. Let all the TT-cores except $X^{(p)}(i_p)$ be ``frozen''.
Since the frame matrix is assumed to be orthogonal, the minimization problem \eqref{eq:rq} is reduced to a smaller problem
\begin{equation}\label{eq:loc}
 \hat X^{(p)}_\new = \arg\min_{\hat X^{(p)}} \trace (  (\hat X^{(p)})^\trans \X_{\neq p}^\trans A \X_{\neq p} \hat X^{(p)} ), \qquad
 \mbox{s.t.}\quad (\hat X^{(p)})^\trans \hat X^{(p)} = I,
\end{equation}
which is in fact equivalent to the standard block eigenproblem for a small matrix.
This~\emph{local} problem is often called~\emph{one--dimensional} (one--site), since it corresponds to the particular mode $p$, and the number of unknowns is linear in the mode size $n_p.$
It is natural to organize the optimization process into \emph{sweeps}, see Alg.~\ref{alg}. In this way, the orthogonalization of the TT-cores is very cheap.

\begin{algorithm}[t]
 \caption{ALS optimization for the block Rayleigh quotient in the TT format} \label{alg}
 \begin{algorithmic}[1]
  \REQUIRE Matrix $A$ and initial guess $X$ in the TT--format~\eqref{eq:att} and~\eqref{eq:btt}, resp.
  \ENSURE Improved eigenvectors $X=[x_b]$ in the TT--format~\eqref{eq:btt}
  \WHILE{stopping criterion is not satisfied}
  \STATE Move index $b$ to $\hat X^{(1)}(i_1,b)$ and make $X^{(2)},\ldots,X^{(d)}$ right--orthogonal
  \FOR[Left--to--right half--sweep]{$p=1,\ldots,d-1$}
    \STATE Solve one--dimensional problem~\eqref{eq:loc} and replace $\hat X^{(p)}:=\hat X^{(p)}_\new$
    \STATE Decompose $\hat X^{(p)}$ by~\eqref{eq:lr} with accuracy $\eps$ and replace $r_p:=r'_p$
    \STATE Ensure the left--orthogonality of $X^{(p)}$
    \STATE In~\eqref{eq:gtt}, merge blocks $G(b)X^{(p+1)}(i_{p+1})=\hat X^{(p+1)}(i_{p+1},b)$
  \ENDFOR
  \STATE Perform right--to--left half--sweep in the same way
  \ENDWHILE
 \end{algorithmic}
\end{algorithm}

The matrix-by-vector product for the local problem scales as $\O(r^2 r_A^2) \matvec(n)$ and $\O(r^3 n r_A)$ additional operations%
\footnote{In complexity estimates we assume that all mode sizes are $\O(n),$ TT--ranks of vector $X$ are $\O(r)$ and TT--ranks of $A$ are $\O(r_A)$}%
, where $\matvec(n_p)$ is a cost of each multiplication of a $n_p \times n_p$ matrix by a vector,
$$
y(i_p) = \sum_{j_p=1}^n A^{(p)}_{\gamma_{p-1},\gamma_p}(i_p,j_p) x(j_p).
$$
For applications in quantum chemistry and physics the mode size is usually not large, and $\O(n^2)$ complexity is acceptable.
In applications of numerical linear algebra and scientific computing $n$ can grow up to thousands, and the use of sparsity, whether possible, is essential to reduce the complexity to linear in $n.$
The $B$ extreme eigenvectors can be found using the block Krylov techniques (e.g. the notable LOBPCG method~\cite{knyazev-lobpcg-2001}) using $\O(B)$ local matvecs and $\O(B^2)$ orthogonalizations of vectors  $x_b.$
The low--rank approximation~\eqref{eq:lr} can be done by the SVD in $\O(Br^3n\min(B,n))$ operations, which can be the slowest part of the algorithm if $n$ and $B$ are both large.
To speed up this step, we can use the cross interpolation of matrices~\cite{tee-cross-2000} (for more robust implementation see~\cite[Alg. 3]{ost-chem-2010}), or incomplete Cholesky method applied to the Gram matrix, see e.g.~\cite{sav-rr-2009}.
Using this methods, we reduce the complexity to $\O(Br^3n),$ that is to the level of other steps of the algorithm.

Summarizing the above, the overall storage and complexity of Alg.~\ref{alg} are
\begin{equation}\label{eq:cost}
 \begin{split}
 \mathrm{mem} & = \O((d+B)nr^2),
 \\
 \work & = \underbrace{\O(dBr^3r_An)+\O(dBr^2r_A^2) \matvec(n)}_{\mbox{\scriptsize local~matvecs}}
         +\underbrace{\O(dB^2 r^2 n)}_{\mbox{\scriptsize orthogonalization}}
         +\underbrace{\O(dB r^3 n)}_{\mbox{\scriptsize truncation}}.
 \end{split}
\end{equation}
Typically the TT--rank of $X$ grows with $B$ at least as $r=\O(B),$ since the TT--format~\eqref{eq:btt} has to represent $B$ vectors with possibly different structures simultaneously.
In this case the first term is more likely to dominate asymptotically for large $B$ and $r$ and we have $\work = \O(dBr^3r_An).$

\begin{table}[t]
 \begin{center}
  \begin{tabular}{cc|cc}
  \multicolumn{2}{c|}{TT notation} & \multicolumn{2}{c}{DMRG literature} \\ \hline
   mode size & $n$ & local dimension & $d$ \\
   dimension & $d$ & number of sites & $n$ \\
   TT--ranks & $r$ & bond dimension  & $D$ \\
   $\#$ eigenvectors & $B$ & $\#$ target states & $M$
  \end{tabular}
 \end{center}
 \caption{Correspondence between notation of our paper and notation used in quantum physics, e.g.~\cite{pizorn-eigb-2012}}
 \label{tab:not}
\end{table}

The complexity of the proposed method can be compared with the one of the algorithms used in quantum physics community, e.g. the DMRG algorithm with~\emph{targeting}~\cite{white-dmrg-target-1993} and the \emph{variational} numerical renormalization group (NRG)~\cite{pizorn-eigb-2012}.
In the notation of this paper (see Table~\ref{tab:not}), the complexity of the DMRG reported in~\cite{pizorn-eigb-2012} is $\O(dBr^3n^3)$ and the complexity of the variational NRG is $\O(Br^3 + dr^3n^3).$
Note that the complexity of the DMRG method is cubic in $n,$ since the optimization over two modes is used in each step.
This is a usual way to couple the subspaces corresponding to different modes (sites) and increase (or decrease) the TT--rank (bond dimension) between the blocks.
In our method the adaptation of the TT--rank $r_k$ is done by the approximate decomposition~\eqref{eq:lr} of the TT--block $\hat X^{(p)}$ which is seen as a $r_{p-1} n_p \times B r_{p+1}$ matrix.

The truncation step introduces the perturbation to the local vector $\hat X^{(p)}.$
The orthogonality of the interface $\X_{\neq p}$ is essential in this step to control the accuracy of the vectors in $X.$
The orthogonality between $x_b$ may be also perturbed, but can be easily recovered by the reorthogonalization of $G(b)=[G_{\alpha'_p,\alpha_p}(b)]$ as a collection of $n$ vectors of size $r'_p r_p.$
In practice this step is not necessary, since the non--orthogonal $G(p)$ is accumulated into $\hat X^{(p+1)},$ which is then replaced by the orthogonal set of the extreme eigenvectors of the local problem in the next step of Alg.~\ref{alg}.

\section{Laplace operator (particle in a box)} \label{sec:lap}
We consider the eigenstates of a particle in a $d$--dimensional box and discretize this problem on a uniform tensor product grid with $n$ elements in each direction as follows
\begin{equation} \label{eq:lap}
-\Delta X = X\Lambda, \qquad \Delta = D \otimes I \otimes \cdots \otimes I + \cdots + I \otimes \cdots \otimes I \otimes D,
\end{equation}
where $I$ is the identity $n\times n$ matrix, and $D = \mathop{\mathrm{tridiag}}(1,-2,1)$ is the standard finite difference discretization of the one-dimensional Laplace operator.
This problem is a nice sanity test for eigensolvers, because the analytical solution of the eigenproblem is available.
The eigenpairs $\{\mu_b,u_b\}$ of the matrix $(-D)$ read
$$
\mu_b = 4 \sin^2\left(\frac{\pi (b+1)}{2(n+1)}\right),
\quad
u_b(j) = \sin\left(\frac{\pi (b+1) (j+1)}{n+1}\right),
\qquad b,j=0,\ldots,n-1,
$$
and for the high--dimensional problem~\eqref{eq:lap} we obtain
\begin{equation}\label{eq:spectrum}
\lambda_{\sigma(b_1, \ldots, b_d)} = \mu_{b_1} + \cdots + \mu_{b_d},
\quad
x_{\sigma(b_1, \ldots, b_d)} = u_{b_1} \otimes \cdots \otimes u_{b_d},
\end{equation}
where the order function $\sigma$ sorts the eigenvalues from low to high.
Alternatively, we can enumerate the~\emph{energy levels} by multiindices as follows:
$$
E_0=\lambda_{\sigma(0,0,\ldots,0)}, \qquad E_1=\lambda_{\sigma(1,0,\ldots,0,0)}=\ldots=\lambda_{\sigma(0,0,\ldots,0,1)},
$$
and so on, and define the invariant eigenspaces $X_k,$ which consist of the eigenvectors corresponding to the same eigenvalues.
The excited eigenstates possess a strong multiplicity, which we show in Table~\ref{tab:laplace}, where $\lambda_0 = d\mu_0$ denotes the ground state energy.
The multiplicity (or \emph{degeneracy}) of the energy levels is a known issue for the convergence of eigensolvers, which makes this example particularly instructive to consider.

\begin{table}[t]
\centering
\begin{tabular}{c|c|c|c}
level $k$& $E_k-E_0$            & multiplicity of $X_k$    & example of eigenvector from $X_k$ \\ \hline
0        & $0$                  & $1$              & $u_0 \otimes u_0 \otimes u_0 \otimes u_0 \otimes \ldots \otimes u_0$ \\
1        & $\mu_1-\mu_0$        & $d$              & $u_1 \otimes u_0 \otimes u_0 \otimes u_0 \otimes \ldots \otimes u_0$ \\
2        & $2(\mu_1-\mu_0)$     & $d(d-1)/2$       & $u_1 \otimes u_1 \otimes u_0 \otimes u_0 \otimes \ldots \otimes u_0$ \\
3        & $\mu_2-\mu_0$        & $d$              & $u_2 \otimes u_0 \otimes u_0 \otimes u_0 \otimes \ldots \otimes u_0$ \\
4        & $3(\mu_1-\mu_0)$     & $d(d-1)(d-2)/6$  & $u_1 \otimes u_1 \otimes u_1 \otimes u_0 \otimes \ldots \otimes u_0$ \\
\end{tabular}
\caption{Low--lying eigenvalues and dimensions of the corresponding invariant subspaces, Laplace example~\eqref{eq:lap}}
\label{tab:laplace}
\end{table}
\begin{figure}[t]
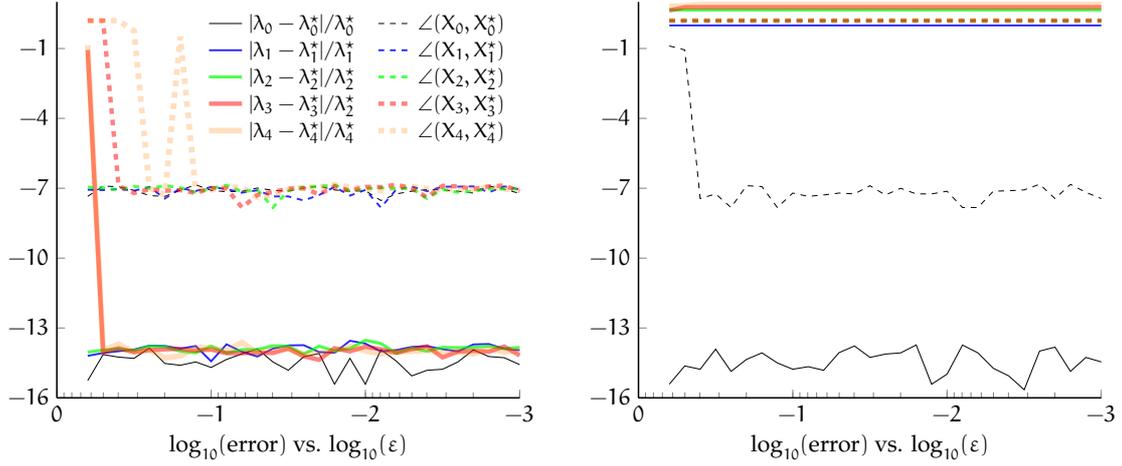

\centering
\resizebox{0.49\textwidth}{!}{\input{./Pic/pgfcol.sty} \begin{tikzpicture}
\begin{axis}[%
  xmode=log,ymode=log,
  cycle list name=eigb,
  xmin=1e-3, xmax=1,
  ymin=1e-16, ymax=10,
  yminorticks=true,
  x dir=reverse,
  legend columns=5,transpose legend,
  legend style={at={(0.99,0.99)},anchor=north east, 
               /tikz/every even column/.append style={column sep=3mm}
               }
  ]

  \pgfplotsset{xlabel/.add={$\log_{10}(\mbox{error})$ vs. $\log_{10}(\eps)$}}

  \addplot+[] table[header=true, x=tol, y=L0]{./dat/laplace_errs_block.dat}; \addlegendentry{$|\lambda_0-\lambda_0^\new|/\lambda_0^\new$};
  \addplot+[] table[header=true, x=tol, y=L1]{./dat/laplace_errs_block.dat}; \addlegendentry{$|\lambda_1-\lambda_1^\new|/\lambda_1^\new$};
  \addplot+[] table[header=true, x=tol, y=L2]{./dat/laplace_errs_block.dat}; \addlegendentry{$|\lambda_2-\lambda_2^\new|/\lambda_2^\new$};
  \addplot+[] table[header=true, x=tol, y=L3]{./dat/laplace_errs_block.dat}; \addlegendentry{$|\lambda_3-\lambda_3^\new|/\lambda_2^\new$};
  \addplot+[] table[header=true, x=tol, y=L4]{./dat/laplace_errs_block.dat}; \addlegendentry{$|\lambda_4-\lambda_4^\new|/\lambda_4^\new$};
  \pgfplotsset{cycle list shift=-5}
  \addplot+[dashed] table[header=true, x=tol, y=A0]{./dat/laplace_errs_block.dat}; \addlegendentry{$\angle(X_0,X_0^\new)$};
  \addplot+[dashed] table[header=true, x=tol, y=A1]{./dat/laplace_errs_block.dat}; \addlegendentry{$\angle(X_1,X_1^\new)$};
  \addplot+[dashed] table[header=true, x=tol, y=A2]{./dat/laplace_errs_block.dat}; \addlegendentry{$\angle(X_2,X_2^\new)$};
  \addplot+[dashed] table[header=true, x=tol, y=A3]{./dat/laplace_errs_block.dat}; \addlegendentry{$\angle(X_3,X_3^\new)$};
  \addplot+[dashed] table[header=true, x=tol, y=A4]{./dat/laplace_errs_block.dat}; \addlegendentry{$\angle(X_4,X_4^\new)$};
    
\end{axis}
\end{tikzpicture}%
} \hfill
\resizebox{0.49\textwidth}{!}{\input{./Pic/pgfcol.sty} \begin{tikzpicture}
\begin{axis}[%
  xmode=log,ymode=log,
  cycle list name=eigb,
  xmin=1e-3, xmax=1,
  ymin=1e-16, ymax=10,
  yminorticks=true,
  x dir=reverse,
  legend columns=5,transpose legend,
  legend style={at={(0.99,0.99)},anchor=north east, 
               /tikz/every even column/.append style={column sep=3mm}
               }
  ]

  \pgfplotsset{xlabel/.add={$\log_{10}(\mbox{error})$ vs. $\log_{10}(\eps)$}}

  \addplot+[] table[header=true, x=tol, y=L0]{./dat/laplace_errs_def.dat};
  \addplot+[] table[header=true, x=tol, y=L1]{./dat/laplace_errs_def.dat};
  \addplot+[] table[header=true, x=tol, y=L2]{./dat/laplace_errs_def.dat};
  \addplot+[] table[header=true, x=tol, y=L3]{./dat/laplace_errs_def.dat};
  \addplot+[] table[header=true, x=tol, y=L4]{./dat/laplace_errs_def.dat};
  \pgfplotsset{cycle list shift=-5}
  \addplot+[dashed] table[header=true, x=tol, y=A0]{./dat/laplace_errs_def.dat};
  \addplot+[dashed] table[header=true, x=tol, y=A1]{./dat/laplace_errs_def.dat};
  \addplot+[dashed] table[header=true, x=tol, y=A2]{./dat/laplace_errs_def.dat};
  \addplot+[dashed] table[header=true, x=tol, y=A3]{./dat/laplace_errs_def.dat};
  \addplot+[dashed] table[header=true, x=tol, y=A4]{./dat/laplace_errs_def.dat};
\end{axis}
\end{tikzpicture}%
}
\caption{Errors in eigenvalues and vectors vs. the thruncation threshold for the Laplace example~\eqref{eq:lap}.
         Left: block-TT Alg.~\ref{alg}. Right: deflation method.
         The reference values $\lambda_k^\new$ and $X_k^\new$ are computed by~\eqref{eq:spectrum}.
         Parameters: dimension $d=5,$ mode size $n=16.$
         }
\label{fig:laplace}
\end{figure}

The matrix $\Delta$ is represented in the TT format with TT--ranks not larger than two, see~\cite{khkaz-lap-2012}.
Each eigenvector $x_{\sigma(b_1,\ldots,b_d)}$ has a rank-one decomposition by~\eqref{eq:spectrum}, and therefore $B$ eigenstates are represented simultaneously in the block--TT format~\eqref{eq:btt} with TT--ranks not larger than $B.$
Since $\{u_b\}$ are orthogonal, the maximal TT--rank grows linearly with $B$ (for small $B$ it is easy to check this directly, but in general one would obtain a complicated combinatorial formula).
Summarizing the above, we see that Alg.~\ref{alg} can be applied to compute lower--lying eigenstates of~\eqref{eq:lap} with $B\lesssim 100.$

We compare the proposed method with the deflation technique, which computes the eigenstates one--by--one using the standard DMRG algorithm.
When $B-1$ eigenvectors $X=[x_b]$ are already computed, the Rayleigh quotient $x^\trans Ax$ is minimized over  the normalized vectors orthogonal to $X.$
All vectors are represented simultaneously in the block--TT format~\eqref{eq:btt}, and the local problem is written as a \emph{deflated} eigenvalue problem with $B-1$ orthogonality constrains.
We could expect the deflation method to be more efficient since the size of the variational problem is smaller.
However, in contrary with the ordinary case, the DMRG combined with the deflation may converge to a wrong eigenpair, corresponding to a larger energy level.

We apply both methods to find $B=30$ low--lying eigenstates of~\eqref{eq:lap} in dimension $d=5.$
This should be enough for the subspaces $X_0,\ldots,X_3,$ and find several vectors from $X_4$.
We vary the truncation parameter $\eps$ and measure the accuracy of the computed eigenvalues comparing them with the analytical values~\eqref{eq:spectrum}.
For $X_k,$ we track the angle between the computed and analytical eigenspace,
$$
\angle(X,Y) = \max \angle(x,y), \qquad x\in \Span X, \quad y \in \Span Y,
$$
using the formula
$
\cos\angle(X,Y) = \sigma_{\min}(X^\trans Y),
$
where $\sigma_{\min}$ is the smallest singular value.

The results are shown in Fig. \ref{fig:laplace}.
We see that all errors in the block method drop down to the machine precision rapidly%
\footnote{Note that the minimal possible error for the angle is the square root of the machine precision due to the $\mathrm{acos}$ function}.
The eigenvectors are represented exactly in the TT format as soon as the rank is large enough, even for very rough truncation threshold $\eps.$
We see that the block DMRG method computes all $B=30$ eigenstates correctly.
The deflation method determines accurately the ground state $X_0$ and two vectors from the first excited subspace $X_1.$
All further eigenvectors belong to the subspaces higher than $X_4,$ and information on the multiplicity and other intermediate states is not computed correctly.
This shows that Alg.~\ref{alg} is reliable in the case of large multiplicities, and therefore can be also applied to compute the eigenstates with small spectral gaps.

\section{Henon--Heiles potential} \label{sec:hh}
We consider a particle in a $d$--dimensional Henon--Heiles potential which is used as a benchmark for high-dimensional quantum molecular dynamics computations \cite{meyer-henon-2002}.
The wavefunction $\psi=\psi(q_1,\ldots,q_d)$ satisfies a stationary Schr\"odinger equation $H\psi = E\psi,$ where the Hamiltonian is defined as follows
\begin{equation}\label{eq:hh}
  \def\Hlap{-\frac12 \Delta}
  \def\Hhar{\frac12 \sum_{k=1}^d q^2_k}
  \def\Hanh{\sum_{k=1}^{d-1}\left(q^2_k q_{k+1} - \frac13 q^3_{k+1} \right)}
    H = \lefteqn{\overbrace{\phantom{\Hlap + \Hhar}}^{\textrm{harmonic part}}}\Hlap + \underbrace{\Hhar + \overbrace{\lambda \Hanh}^{\textrm{anharmonic part}}}_{\textrm{Henon-Heiles potential}~V(q_1,\ldots,q_d)},
\end{equation}
where parameter $\lambda=0.111803$ controls the anharmonic contribution to the potential.

The principal part of the Henon-Heiles operator describes a harmonic oscillator, whose eigenstates are products of a Gaussian by Hermite polynomials and have tensor rank $1$.
Therefore, for moderate anharmonicity $\lambda$ the rank-$1$ basis of eigenfunctions of the harmonic oscillator is a good choice for the discretization of \eqref{eq:hh}.
The Galerkin discretization scheme results not only in dense stiffness and mass matrices, but also in a dense matrix that describes the action of the potential $V$.
The DVR(discrete variable representation) scheme uses a collocation of $V$ and $\psi$ on the nodes of the Hermite polynomials and is known to provide the same order of accuracy for the eigenproblem as the Galerkin method.
\begin{figure}[t]
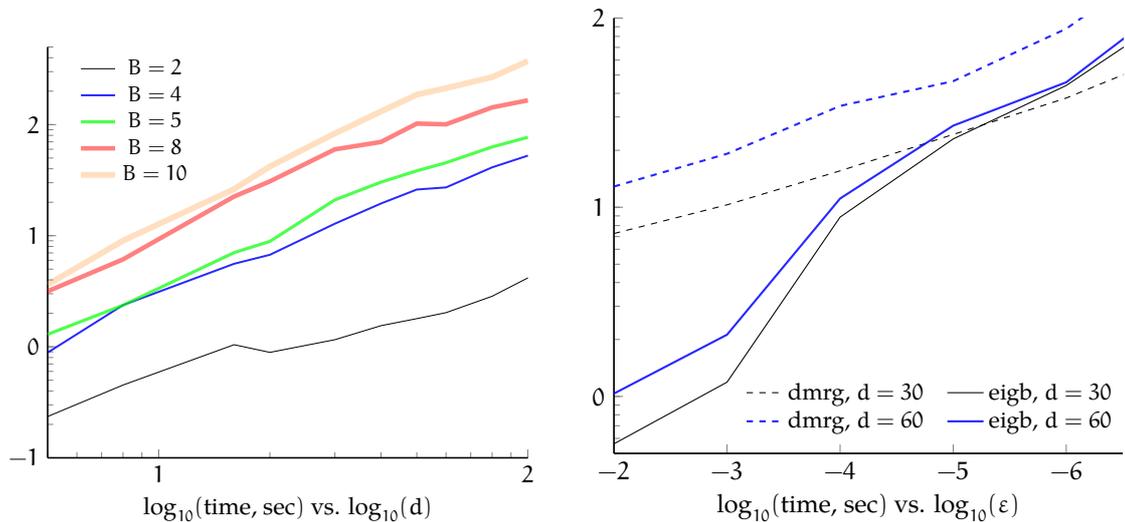

\centering
\resizebox{0.49\textwidth}{!}{\input{./Pic/pgfcol.sty} \begin{tikzpicture}
\begin{axis}[%
  xmode=log,ymode=log,
  cycle list name=eigb,
  xmin=0, xmax=100,
  ymin=0.1, ymax=5e2,
  yminorticks=true,
  legend style={at={(0.05,0.99)},anchor=north west}]

  \pgfplotsset{xlabel/.add={$\log_{10}(\mbox{time, sec})$ vs. $\log_{10}(d)$}}
  \addplot+[] table[header=false, x index=0, y index=1]{./dat/hh_ttimes_L_B_e3.dat}; \addlegendentry{$B=2$};
  \addplot+[] table[header=false, x index=0, y index=2]{./dat/hh_ttimes_L_B_e3.dat}; \addlegendentry{$B=4$};
  \addplot+[] table[header=false, x index=0, y index=3]{./dat/hh_ttimes_L_B_e3.dat}; \addlegendentry{$B=5$};
  \addplot+[] table[header=false, x index=0, y index=4]{./dat/hh_ttimes_L_B_e3.dat}; \addlegendentry{$B=8$};
  \addplot+[] table[header=false, x index=0, y index=5]{./dat/hh_ttimes_L_B_e3.dat}; \addlegendentry{$B=10$};
\end{axis}
\end{tikzpicture}%
} \hfill
\resizebox{0.49\textwidth}{!}{\input{./Pic/pgfcol.sty} \begin{tikzpicture}
\begin{axis}[%
  xmode=normal,ymode=log,
  cycle list name=eigb,
  xtick={-6,...,-2},
  xmin=-6.5,xmax=-2,
  ymin=0.5, ymax=1e2,
  yminorticks=true,
  x dir=reverse,
  legend columns=2,transpose legend,
  legend style={at={(1,0.03)},anchor=south east,
               /tikz/every even column/.append style={column sep=3mm}
               }
  ]

  \pgfplotsset{xlabel/.add={$\log_{10}(\mbox{time, sec})$ vs. $\log_{10}(\eps)$}}
  \addplot+[dashed] table[header=true, x index = 0, y index = 2]{./dat/hh_ttimes2_e_L.dat}; \addlegendentry{dmrg, $d=30$};
  \addplot+[dashed] table[header=true, x index = 0, y index = 4]{./dat/hh_ttimes2_e_L.dat}; \addlegendentry{dmrg, $d=60$};
  \pgfplotsset{cycle list shift=-2}
  \addplot+[] table[header=false, x index=0, y index=1]{./dat/hh_ttimes_e_L_B1.dat}; \addlegendentry{eigb, $d=30$};
  \addplot+[] table[header=false, x index=0, y index=2]{./dat/hh_ttimes_e_L_B1.dat}; \addlegendentry{eigb, $d=60$};
\end{axis}
\end{tikzpicture}%
}
\caption{CPU time vs. dimension $d$ (left) and truncation threshold $\eps$ (right) for the Henon-Heiles example~\eqref{eq:hh}.
         Parameters: $\eps=10^{-3}$ (left), $B=2$ (right), $n=28$.}
\label{fig:hh_ttimes}
\end{figure}

\begin{figure}[t]
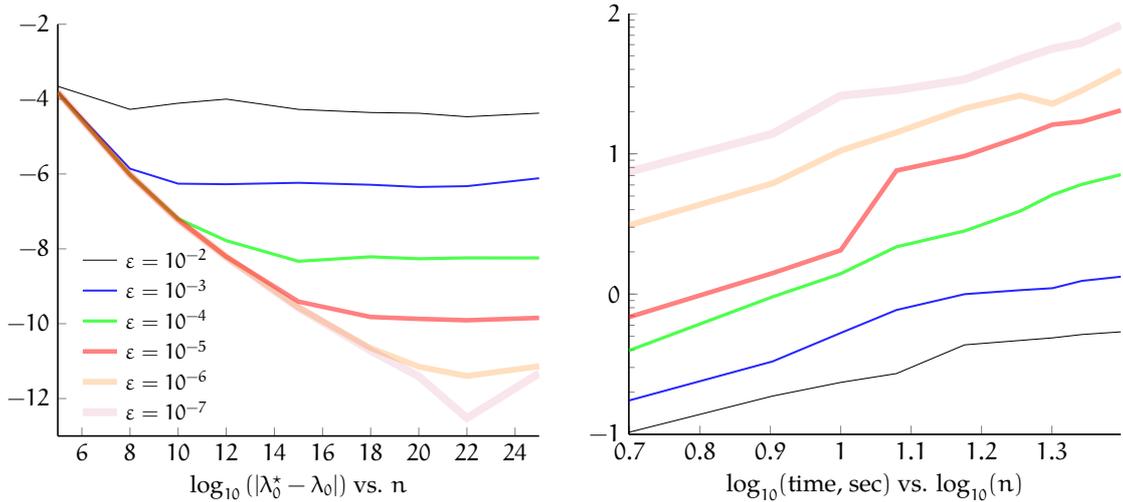

\centering
\resizebox{0.49\textwidth}{!}{\input{./Pic/pgfcol.sty} \begin{tikzpicture}
\begin{axis}[%
  xmode=normal,ymode=log,
  cycle list name=eigb,
  xmin=5, xmax=25,
  ymin=1e-13, ymax=1e-2,
  yminorticks=true,
  legend style={at={(0.03,0.01)},anchor=south west}]

  \pgfplotsset{xlabel/.add={$\log_{10}\left(\left|\lambda^\new_0-\lambda_0\right|\right)$ vs. $n$}}
  \addplot+[] table[header=false, x index=0, y index=1]{./dat/hh_errors_L30_e_n.dat}; \addlegendentry{$\eps=10^{-2}$};
  \addplot+[] table[header=false, x index=0, y index=2]{./dat/hh_errors_L30_e_n.dat}; \addlegendentry{$\eps=10^{-3}$};
  \addplot+[] table[header=false, x index=0, y index=3]{./dat/hh_errors_L30_e_n.dat}; \addlegendentry{$\eps=10^{-4}$};
  \addplot+[] table[header=false, x index=0, y index=4]{./dat/hh_errors_L30_e_n.dat}; \addlegendentry{$\eps=10^{-5}$};
  \addplot+[] table[header=false, x index=0, y index=5]{./dat/hh_errors_L30_e_n.dat}; \addlegendentry{$\eps=10^{-6}$};
  \addplot+[] table[header=false, x index=0, y index=6]{./dat/hh_errors_L30_e_n.dat}; \addlegendentry{$\eps=10^{-7}$};
\end{axis}
\end{tikzpicture}%
}\hfill
\resizebox{0.49\textwidth}{!}{\input{./Pic/pgfcol.sty} \begin{tikzpicture}
\begin{axis}[%
  xmode=log,ymode=log,
  cycle list name=eigb,
  xmin=5, xmax=25,
  ymin=0.1, ymax=1e2,
  yminorticks=true,
  legend style={at={(0.05,0.99)},anchor=north west}]

  \pgfplotsset{xlabel/.add={$\log_{10}(\mbox{time, sec})$ vs. $\log_{10}(n)$}}
  \addplot+[] table[header=false, x index=0, y index=1]{./dat/hh_ttimes_L30_e_n.dat}; 
  \addplot+[] table[header=false, x index=0, y index=2]{./dat/hh_ttimes_L30_e_n.dat}; 
  \addplot+[] table[header=false, x index=0, y index=3]{./dat/hh_ttimes_L30_e_n.dat}; 
  \addplot+[] table[header=false, x index=0, y index=4]{./dat/hh_ttimes_L30_e_n.dat}; 
  \addplot+[] table[header=false, x index=0, y index=5]{./dat/hh_ttimes_L30_e_n.dat}; 
  \addplot+[] table[header=false, x index=0, y index=6]{./dat/hh_ttimes_L30_e_n.dat}; 
\end{axis}
\end{tikzpicture}%
}
\caption{Eigenvalue errors (left) and CPU times (right) vs. mode size $n$ and truncation threshold $\eps$ for the Henon-Heiles example~\eqref{eq:hh}.
         The reference value $\lambda^\new$ is computed by the same algorithm with $\eps=10^{-7}, n=28.$
         Parameters: $B=2,$ $d=30$}
\label{fig:hh_n}
\end{figure}

The one-dimensional Laplace operator $D = -\frac{d^2}{dq^2}$ in the Hermite-DVR approach is discretized (see~\cite{baye-dvr-1986}) as follows
\begin{equation}\nonumber
  D_{ij} = \begin{cases}
            \frac{1}{6}(4 n - 1  - 2t^2_i),               & \quad i = j, \\
           (-1)^{(i-j)}(2(t_i - t_j)^{-2} - \frac{1}{2}), & \quad i \neq j,
          \end{cases}
\end{equation}
where $t_1,\ldots,t_n$ are the roots of the $n$-th Hermite polynomial, i.e. the \emph{Hermite mesh}.
The discretization of the $d$-dimensional Laplace operator is written in the same way as in~\eqref{eq:lap} and has TT--ranks not larger than two.
The potential $V$ is constructed as a sum of rank-1 monomials and can be represented by the TT format with TT--ranks not larger than $7,$ see \cite{khos-dmrg-2010}.
Therefore, Alg.~\ref{alg} can be applied to find the low--lying eigenstates of~\eqref{eq:hh}.

As an initial guess, we take a random TT tensor of rank $B$ in the block form \eqref{eq:btt}.
We test the block TT algorithm for different dimensions $d$, Frobenius-norm truncation tolerances for eigenvalues $\eps$ and number of eigenstates $B$.
Numerical experiments show that the internal convergence of the eigenvalues is quadratic w.r.t. the truncation threshold $\eps,$ as expected from the perturbation theory.
The errors in the eigenvalues do not grow with the dimension, so the method can be exploited for higher--dimensional systems.
From Fig.~\ref{fig:hh_ttimes} (left), we observe that the cost grows mildly with the dimension and $B$, which allows to compute energy levels with high accuracy.

In Fig.~\ref{fig:hh_ttimes} (right) we compare our block-TT method (Alg.~\ref{alg}, referenced as `eigb') with the DMRG technique from~\cite{khos-dmrg-2010}, intended for searching one lowest eigenpair only.
Alg.~\ref{alg} computes at least $B=2$ eigenstates, and the TT-ranks of the block--TT format which represents two eigenstates are larger (sometimes significantly) than the TT-ranks of the ground state returned by the DMRG method.
We see how large ranks of the excited state manifest themselves for high accuracies in Fig.~\ref{fig:hh_ttimes}.
Despite the lower cost w.r.t. the mode size $n$, the `eigb' method requires the same CPU time as the DMRG due to larger TT-ranks, which leads to a higher complexity of the local eigenproblem~\eqref{eq:loc}.

Finally we investigate the performance of Alg.~\ref{alg} w.r.t. the mode size (number of Hermite polynomials) $n$, see Fig. \ref{fig:hh_n}.
The error decreases exponentially with $n$ due to the Hermite-DVR scheme, but only until the level $\O(\eps^2)$ governed by the tensor truncation threshold.
The CPU time growth rate is balanced between $\O(n)$ (truncation step) and $\O(n^2)$ (dense matrix-by-vector multiplication).

\section{Heisenberg model} \label{sec:spin}
The one-dimensional Heisenberg model is one of the classical applications of the MPS/DMRG algorithms, so it is worth to test our approximate block eigenvalue solver on it, and compare it with the established software from the MPS community.
This model describes the interaction of spins on a one-dimensional lattice.
The Hamiltonian for the antiferromagnetic case is written in the following form
\begin{equation}\label{eq:heisen}
 H = \sum_{i=1}^{d-1} \S_i \S_{i+1},
\end{equation}
where $\S_i$ is the spin operator.
The product $\S_i \S_{i+1}$ can be written in terms of the spin components as follows
\begin{equation}\label{eq:spin1}
\S_i \S_{i+1} = (S_x)_i (S_x)_{i+1} + (S_y)_i (S_y)_{i+1} + (S_z)_i (S_z)_{i+1}.
\end{equation}
Since the operator $(S_x)_i$ affects only the $i$-th spin, it has the form
\begin{equation}\label{eq:spin2}
   (S_x)_i = I \otimes \ldots \otimes S_x \otimes \ldots \otimes I,
\end{equation}
and similarly for $(S_y)_i$, $(S_z)_i$.
The ``elementary'' spin operators $S_x$, $S_y$ and $S_z$ are the Pauli matrices, defined as follows
\begin{equation}\label{eq:spin3}
  S_x = \frac12 \begin{pmatrix} 0 & 1 \\  1 & 0 \\ \end{pmatrix}, \quad
  S_y = \frac12 \begin{pmatrix} 0 & -i \\ i & 0 \\ \end{pmatrix}, \quad
  S_z = \frac12 \begin{pmatrix} 1 & 0 \\ 0 & -1 \\ \end{pmatrix}.
\end{equation}
Using the equations \eqref{eq:heisen},\eqref{eq:spin1},\eqref{eq:spin2} and \eqref{eq:spin3},
it is very easy to construct the TT representation of the Heisenberg Hamiltonian.
Moreover, it is also not difficult to show that its TT ranks are not larger than $5,$ see~\cite{dkh-cme-2012}.

\begin{figure}[t]
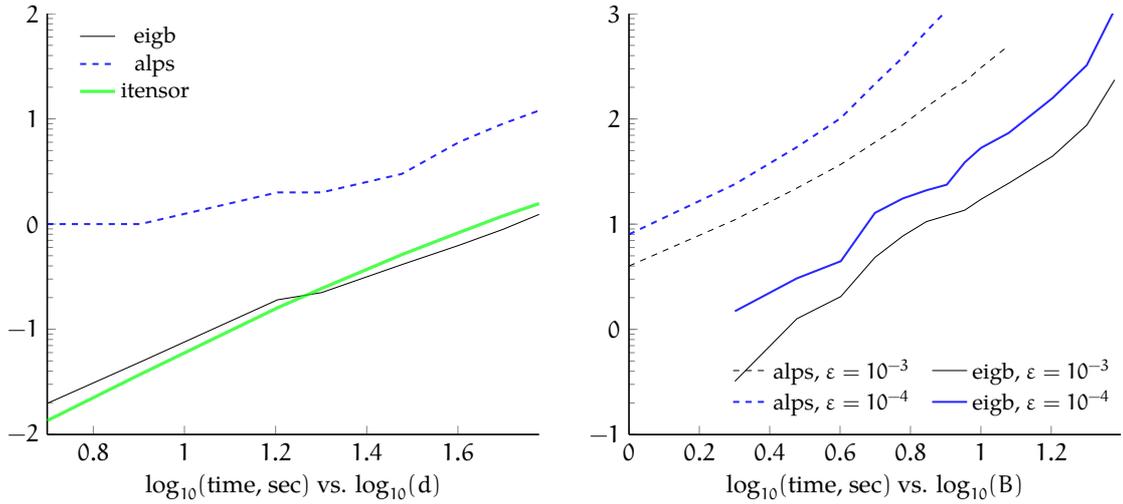

\centering
\resizebox{0.49\textwidth}{!}{\input{./Pic/pgfcol.sty} \begin{tikzpicture}
\begin{axis}[%
  xmode=log,ymode=log,
  cycle list name=eigb,
  xmin=0, xmax=60,
  ymin=0.01, ymax=1e2,
  yminorticks=true,
  legend style={at={(0.05,0.99)},anchor=north west}]

  \pgfplotsset{xlabel/.add={$\log_{10}(\mbox{time, sec})$ vs. $\log_{10}(d)$}}
  \addplot+[] table[header=false, x index=0, y index=1]{./dat/heisen_ttimes_L_m_B1_e3.dat}; \addlegendentry{eigb};
  \addplot+[dashed] table[header=false, x index=0, y index=2]{./dat/heisen_ttimes_L_m_B1_e3.dat}; \addlegendentry{alps};
  \addplot+[] table[header=false, x index=0, y index=3]{./dat/heisen_ttimes_L_m_B1_e3.dat}; \addlegendentry{itensor};
\end{axis}
\end{tikzpicture}%
} \hfill
\resizebox{0.49\textwidth}{!}{\input{./Pic/pgfcol.sty} \begin{tikzpicture}
\begin{axis}[%
  xmode=log,ymode=log,
  cycle list name=eigb,
  xmin=1, xmax=25,
  ymin=0.1, ymax=1e3,
  yminorticks=true,
  legend columns=2,transpose legend,
  legend style={at={(1, 0.03)},anchor=south east,
               /tikz/every even column/.append style={column sep=3mm}
               }
  ]

  \pgfplotsset{xlabel/.add={$\log_{10}(\mbox{time, sec})$ vs. $\log_{10}(B)$}}
  \addplot+[dashed] table[header=false, x index=0, y index=3]{./dat/heisen_ttimes_L30_m_B_e.dat}; \addlegendentry{alps, $\eps=10^{-3}$};
  \addplot+[dashed] table[header=false, x index=0, y index=4]{./dat/heisen_ttimes_L30_m_B_e.dat}; \addlegendentry{alps, $\eps=10^{-4}$};
  \pgfplotsset{cycle list shift=-2}
  \addplot+[] table[header=false, x index=0, y index=1]{./dat/heisen_ttimes_L30_m_B_e.dat}; \addlegendentry{eigb, $\eps=10^{-3}$};
  \addplot+[] table[header=false, x index=0, y index=2]{./dat/heisen_ttimes_L30_m_B_e.dat}; \addlegendentry{eigb, $\eps=10^{-4}$};
\end{axis}
\end{tikzpicture}%
}
\caption{CPU times of different methods vs. $d$ (left), $B$ (right), Heisenberg example.
         Parameters: $\eps=10^{-3}$, $B=1$ (left), $d=30$ (right).
         }
\label{fig:heisen_ttimes}
\end{figure}

\begin{figure}[t]
\centering
\resizebox{0.49\textwidth}{!}{\input{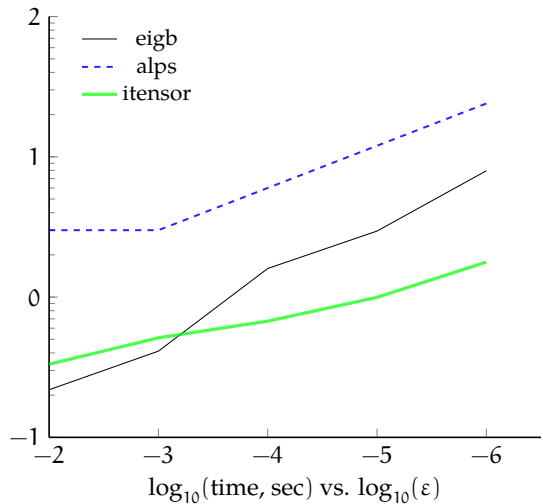} \begin{tikzpicture}
\begin{axis}[%
  xmode=normal,ymode=log,
  cycle list name=eigb,
  xtick={-6,...,-2},
  xmin=-6.5, xmax=-2,
  ymin=0.1, ymax=1e2,
  yminorticks=true,
  x dir=reverse,
  legend style={at={(0.05,0.99)},anchor=north west}]

  \pgfplotsset{xlabel/.add={$\log_{10}(\mbox{time, sec})$ vs. $\log_{10}(\eps)$}}
  \addplot+[] table[header=false, x index=0, y index=1]{./dat/heisen_ttimes_L30_m_B1_e.dat}; \addlegendentry{eigb};
  \addplot+[dashed] table[header=false, x index=0, y index=2]{./dat/heisen_ttimes_L30_m_B1_e.dat}; \addlegendentry{alps};
  \addplot+[] table[header=false, x index=0, y index=3]{./dat/heisen_ttimes_L30_m_B1_e.dat}; \addlegendentry{itensor};
\end{axis}
\end{tikzpicture}%
}
\caption{CPU times of different methods vs. $\eps$, Heisenberg example.
         Parameters: $B=1$, $d=30$.
         }
\label{fig:heisen_ttimes_e}
\end{figure}

\begin{figure}[t]
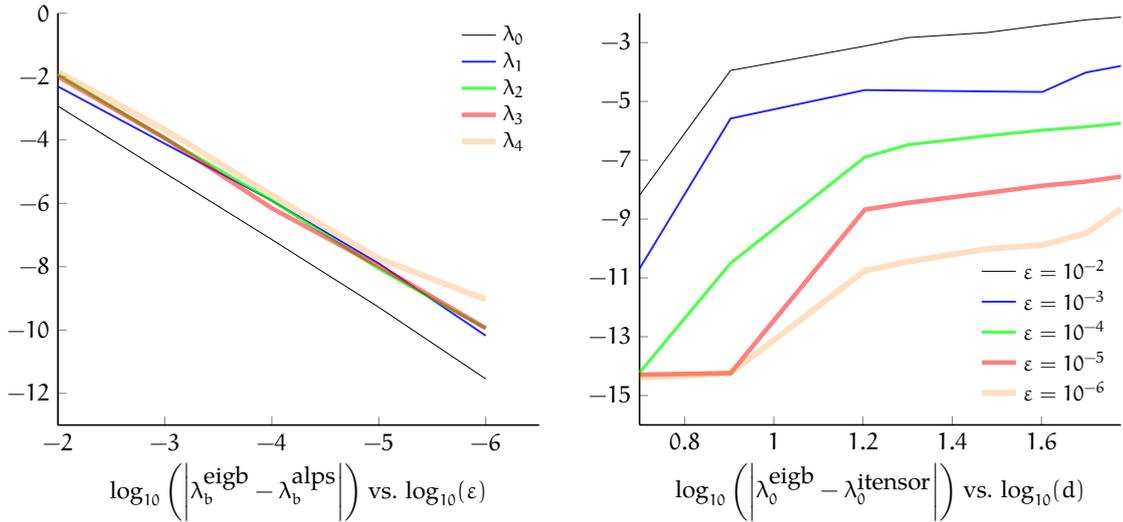

\centering
\resizebox{0.49\textwidth}{!}{\input{./Pic/pgfcol.sty} \begin{tikzpicture}
\begin{axis}[%
  xmode=normal,ymode=log,
  cycle list name=eigb,
  xtick={-6,...,-2},
  xmin=-6.5, xmax=-2,
  ymin=1e-13, ymax=1,
  yminorticks=true,
  x dir=reverse,
  legend style={at={(0.99,0.99)},anchor=north east}]

  \pgfplotsset{xlabel/.add={$\log_{10}\left(\left|\lambda^{\mbox{eigb}}_b-\lambda^{\mbox{alps}}_b\right|\right)$ vs. $\log_{10}(\eps)$}}
  \addplot+[] table[header=false, x index=0, y index=1]{./dat/heisen_errs_L30.dat}; \addlegendentry{$\lambda_0$};
  \addplot+[] table[header=false, x index=0, y index=2]{./dat/heisen_errs_L30.dat}; \addlegendentry{$\lambda_1$};
  \addplot+[] table[header=false, x index=0, y index=3]{./dat/heisen_errs_L30.dat}; \addlegendentry{$\lambda_2$};
  \addplot+[] table[header=false, x index=0, y index=4]{./dat/heisen_errs_L30.dat}; \addlegendentry{$\lambda_3$};
  \addplot+[] table[header=false, x index=0, y index=5]{./dat/heisen_errs_L30.dat}; \addlegendentry{$\lambda_4$};
\end{axis}
\end{tikzpicture}%
} \hfill
\resizebox{0.49\textwidth}{!}{\input{./Pic/pgfcol.sty} \begin{tikzpicture}
\begin{axis}[%
  xmode=log,ymode=log,
  cycle list name=eigb,
  xmin=0, xmax=60,
  ymin=1e-16, ymax=1e-2,
  yminorticks=true,
  legend style={at={(0.99,0.03)},anchor=south east}]

  \pgfplotsset{xlabel/.add={$\log_{10}\left(\left|\lambda^{\mbox{eigb}}_0-\lambda^{\mbox{itensor}}_0\right|\right)$ vs. $\log_{10}(d)$}}
  \addplot+[] table[header=false, x index=0, y index=1]{./dat/heisen_errs_B1.dat}; \addlegendentry{$\eps=10^{-2}$};
  \addplot+[] table[header=false, x index=0, y index=2]{./dat/heisen_errs_B1.dat}; \addlegendentry{$\eps=10^{-3}$};
  \addplot+[] table[header=false, x index=0, y index=3]{./dat/heisen_errs_B1.dat}; \addlegendentry{$\eps=10^{-4}$};
  \addplot+[] table[header=false, x index=0, y index=4]{./dat/heisen_errs_B1.dat}; \addlegendentry{$\eps=10^{-5}$};
  \addplot+[] table[header=false, x index=0, y index=5]{./dat/heisen_errs_B1.dat}; \addlegendentry{$\eps=10^{-6}$};
\end{axis}
\end{tikzpicture}%
}
\caption{Discrepancies between eigenvalues in different methods vs. eigenvalue number (left) and $d$ (right), Heisenberg example.
         Parameters: $d=30,B=5$ (left), $B=1$ (right).
         }
\label{fig:heisen_errors}
\end{figure}

For different lengths of the spin chain we compute the lowest eigenvalues
with the help of our block eigenvalue solver, and compare it with two packages containing the DMRG:
ALPS (Algorithms and Libraries for Physics Simulations)\footnote{\url{http://alps.comp-phys.org/}, Release 2.1.1} and ITensor\footnote{\url{http://itensor.org/}, downloaded on May 20, 2013 (no particular release)}.
The ALPS allows to compute several excited states, whereas the ITensor is devised for targeting the ground state only.
Using Alg.~\ref{alg} we can compute $B\geq2$ eigenstates, so the comparison with ITensor and ALPS with $B=1$ assumes that our algorithm is applied with $B=2$ and the excited state is computed but not used.
To introduce a fair amount of optimization, all software was compiled using Intel C/Fortran compiler 2013 and linked with the optimized \textsc{Lapack}/\textsc{Blas} packages provided in the MKL library.
The experiments have been computed using one Intel Xeon processor with $4$ cores at $2$GHz.

The computational times of Alg.~\ref{alg} (`eigb'), ALPS and ITensor are presented in Figs~\ref{fig:heisen_ttimes} and~\ref{fig:heisen_ttimes_e}.
We observe that all methods manifest the polynomial complexity scaling in the number of spins $d$ and targeting eigenstates $B$.
The dependence on $B$ is of particular interest.
Recalling the complexity estimate~\eqref{eq:cost}, we may expect  $\work=\O(B^4),$ if  the TT--ranks behave as $r \sim B.$
In Fig.~\ref{fig:heisen_ttimes} we see even milder growth of the CPU time, $\time=\O(B^\beta),$ $\beta \approx 2.5,$ and the same phenomenon is observed for the DMRG method from the ALPS.
It can be explained by the non-uniform distribution of the TT--ranks, which makes the complexity estimate~\eqref{eq:cost} larger than the actual computation time.

Considering the memory usage, we observe that the ALPS package failed to compute more than $12$ targets due to the Out-of-Memory exception, while our implementation of Alg.~\ref{alg} goes readily beyond $20$ eigenstates.

The accuracy parameter $\eps$ requires an additional comment, since it is used differently in different methods.
The block TT method performs the Frobenius-norm $\eps$-approximation of the eigenvectors, whereas the other two packages are designed to keep the accuracy of the eigenvalue, and hence are parametrized with $\eps^2$ as a threshold.

We see that the performances of our block TT technique and the method from ITensor are comparable, and the algorithm from ALPS is significantly slower.
Taking into account that our approach allows to compute several eigenstates (which seems to be not the case for ITensor), we may conclude that it overcomes the current state-of-the-art software.
To be sure that the correctness is maintained, we show the difference between the eigenvalues computed by the eigb and the other two methods in Fig. \ref{fig:heisen_errors}, which demonstrates that the $\O(\eps^2)$ accuracy of the eigenvalues is achieved.

\section{Conclusions}
We propose the efficient one--site rank--adaptive algorithm for the computation of several extreme eigenvalues of a high--dimensional Hermitian matrix.
The eigenstates are represented simultaneously by the block tensor train format, which is flexible and has a particular form in each variational step of the algorithm.
The complexity of the proposed method is linear with respect to the mode size for sparse matrices and quadratic for dense matrices, which is better than the complexity of the DMRG algorithm.
The proposed method has the same asymptotic in the number of eigenstates, the TT--ranks and the dimension, as the DMRG and NRG algorithms used in quantum physics.

The algorithm is verified on a number of representative problems.
First, it is tested on the high--dimensional Laplace matrix, which describes the eigenstates of the particle in a box, and it is demonstrated that the algorithm is capable of finding the eigenstates with high multiplicity.
Second, it is applied to the stationary Schr\"odinger equation with the Henon-Heiles potential, and it is shown that the proposed method can compute the ground state faster than the DMRG algorithm.
Finally, we apply the algorithm to the Heisenberg spin chain and show that it is competitive with the state of the art publicly available DMRG implementations.

The framework developed in this paper can be applied to a wider class of problems, which are formulated via the subspace optimization.

\section*{Software implementation}

We have implemented the main Algorithm \ref{alg} in Fortran with interfaces to Python and MATLAB.
The codes are available online at
\begin{itemize}
 \item \url{http://github.com/oseledets/tt-fort}
 \item \url{http://github.com/oseledets/ttpy}
 \item \url{http://github.com/oseledets/TT-Toolbox}
\end{itemize}



\end{document}